\documentclass[12pt,a4paper]{amsart}
\usepackage[latin1]{inputenc}
\usepackage{latexsym}
\usepackage{color,graphicx,shortvrb}
\newtheorem{theorem}{Theorem}

\newtheorem{corollary}[theorem]{Corollary}
\newtheorem{example}{Example}
\title{Metrics with good Corona properties}

\author{J. M. Almira, A. J. L\'{o}pez-Moreno and N. Del Toro
}

\begin{document}

\begin{abstract}
In this note, we give several definitions of metric corona
properties which could be of interest in Set Topology, Functional
Analysis and Approximation Theory, and prove that there are
complete metrizable t.v.s. which are nice in the sense that they
have a metric which is invariant by translations, but they do not
have good corona properties. All classical function spaces satisfy
good corona properties but it is an open question to know if this
also holds for the more general setting of locally convex t.v.s.
\end{abstract}

\maketitle

\section{Introduction. Definition of corona properties}
Let $(\mathbf{X,d})$\textbf{\ }be a metric space. The study of
topologies on the sets
\[
\mathbf{CL(X)=\{A\subset X:A}\text{ is closed and }\mathbf{A}\neq
\emptyset \}
\]
and
\[
\mathbf{CLB(X)=\{A\subset X:A}\text{ is closed, bounded and
}\mathbf{A}\neq \emptyset \}
\]
is of interest in several branches of Mathematics, as Set
Topology, Functional Analysis, etc. (see \cite{beer}). These
topologies are usually defined as weak topologies associated to a
family of functionals which depend on the metric $\mathbf{d}$ in
some sense and, of course, it is of interest to know if certain
natural maps related to the metric $\mathbf{d}$ are or are not
continuous. The corona properties of $(\mathbf{X,d})$, appear in a
natural form for the study of the continuity of the map $\mathcal{B}:\mathbf{X}%
\times \lbrack 0,\infty )\rightarrow \mathbf{CLB(X)}$ defined by $\mathcal{B}%
(\mathbf{x},t)=\overline{\mathbf{B}}(\mathbf{x},t)$, where
$\mathbf{B}(\mathbf{x},t)$ stands for the ball centered at
$\mathbf{x}$ and with radius $t$ for the considered metric
$\mathbf{d}$. It seems that such a map must be
continuous for any reasonable topology on $\mathbf{CL(X)}$ and any metric $%
\mathbf{d}$ but this is not true. For example, if we consider on $\mathbf{%
CLB(X)}$ the classical Hausdorff metric.
\[
\mathbf{H(A,B)}=\max \{\sup_{\mathbf{x\in A}}\mathbf{d(x,B)},\sup_{\mathbf{%
x\in B}}\mathbf{d(x,A)}\}
\]
then the following result can be proved.

\begin{theorem}
Let $\mathbf{x}\in\mathbf{X}$ be fixed, and let us consider the
map $\mathcal{B}_{\mathbf{x}}(t)=\mathcal{B}(\mathbf{x},t)$. Then
$\mathcal{B}_{\mathbf{x}}$ is right-continuous if and only if the
following property holds:

\noindent $\mathbf{(WCP)}$ If
$\lim_{n\rightarrow\infty}\mathbf{d}(\mathbf{x},\mathbf{y}_n)=t$
then
$\lim_{n\rightarrow\infty}\mathbf{d}(\mathbf{y}_n,\overline{\mathbf{B}}_d(\mathbf{x},t))=0$
\end{theorem}

%TCIMACRO{
%\TeXButton{margen}{\noindent%
%}}%
%BeginExpansion
\noindent%
%
%EndExpansion
\textbf{Proof. }Let us assume that $\mathcal{B}_x$ is
right-continuous, and let $(\varepsilon
_{n})\searrow 0^{+}$. Suppose that $\mathbf{y}_{n}\in \mathbf{X}$ satisfies $%
\mathbf{d(x},\mathbf{y}_{n})\leq t_{0}+\varepsilon _{n}$ for all
$n$. It follows from our hypothesis on the continuity of
$\mathcal{B}_{\mathbf{x}}$ and the definition of $\mathbf{H}$ that
\[
\lim_{n\rightarrow \infty }\mathbf{d}\left(\mathbf{y}_{n}\mathbf{\mathbf{,\mathbf{%
\overline{\mathbf{B}}(x}}},t_{0} )\right)\leq \lim_{n\rightarrow
\infty }\mathbf{H}\left(
\overline{\mathbf{B}}(\mathbf{x},t_{0}+\varepsilon _{n}
),\overline{\mathbf{B}}(\mathbf{x} ,t_{0})\right)=0\text{. }
\]
This proves the first implication. Let us now assume that the
property  $(\mathbf{WCP})$ holds. If $\mathcal{B}_{\mathbf{x}}$ is
not right-continuous at $t_0$, then there exists $\varepsilon
>0$ and an infinity sequence of positive real numbers,
$\{t_n\}_{n=0}^{\infty}$ such that
$\lim_{n\rightarrow\infty}t_n=t_0$, $t_n\geq t_0$ for all $n$ and
\begin{eqnarray*}
\varepsilon &\leq&
\mathbf{H}(\mathbf{B}(\mathbf{x_0},t_0),\mathbf{B}(\mathbf{x}_n,t_n))\\
&=&
\sup_{\mathbf{z}\in\mathbf{B}(\mathbf{x_0},t_n)}\mathbf{d(z},\mathbf{B(x}_0,t_0)),\
n\in\mathbb{N},
\end{eqnarray*}
It follows that there exists a sequence
$\{\mathbf{z}_n\}_{n=0}^{\infty}\subset \mathbf{X}$ such that
$\mathbf{d}(\mathbf{x}_0,\mathbf{z}_n)\leq t_0+(t_n-t_0)$, for all
$n$, but $\lim\mathbf{d(z}_n,\mathbf{B(x}_0,t_0))\neq 0$, which is
in contradiction with $\mathbf{(WCP)}$.  $\Box$

We say that the metric space $\mathbf{(X,d)}$ has the weak corona
property if it satisfies the property \textbf{(WCP)} above. Of
course, this property is satisfied by all normed linear
spaces, but it is not satisfied by all metric spaces. For example, $\mathbf{X_1}=%
\mathbb{R}^{2}\setminus \{\mathbf{x=}(x,y):\left\| \mathbf{x}\right\|
\leq 1$ and $x>0\}$ with the usual metric of $\mathbb{R}^{2}$ has not
the weak corona property (see the figure).

%\begin{picture}(200,200)(-120,-5)
\begin{center}
    \includegraphics[scale=0.8]{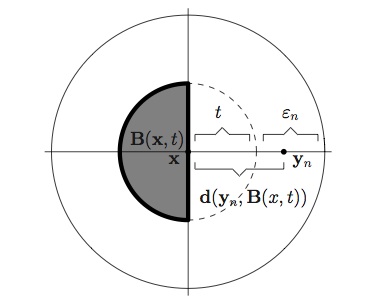}
    \end{center}
   % \put(-32,115){$\varepsilon_n$}
    %\put(-75,115){$t$}
    %\put(-85,60){$\mathbf{d}(\mathbf{y}_n,\mathbf{B}(x,t))$}
    %\put(-130,97){$\mathbf{B}(\mathbf{x},t)$}
    %\put(-105,85){$\mathbf{x}$}
    %\put(-25,85){$\mathbf{y}_n$}
%\end{picture}

%\vskipe

%\begin{center}
%\includegraphics[scale=0.8]{dibujo2.eps}
%\end{center}
Moreover, if we set
$\mathbf{X}_2=\overline{\mathbf{X_1}}^{\mathbf{R}^2}$, then
$\mathbf{X}_2$ has the weak corona property and the map
$\mathcal{B}_0:[0,\infty)\rightarrow\mathbf{CLB}(\mathbf{X}_2)$
given by $\mathcal{B}_0(t)=\overline{\mathbf{B}(0,t)}$ is not
left-continuous.

There are more possibilities to give definitions of corona
properties, and several other results can be proved in connection
with the continuity of the map $\mathcal{B}$ and these properties.
More precisely, we say that $\mathbf{(X},\mathbf{d})$ has the
corona property if there exists a constant $C>0$ such that the relation $%
\mathbf{d}(\mathbf{x},\mathbf{y})\leq t+\varepsilon $ implies the relation
$\mathbf{d}(\mathbf{y},\overline{\mathbf{B}}%
(\mathbf{x},t))\leq C\varepsilon $ for all $t>0$ and all
$\varepsilon \in [0,\varepsilon _{0}(t))$, for a certain
$\varepsilon _{0}(t)>0$. Finally, we say that
$(\mathbf{X},\mathbf{d})$ has the strong corona property if there
exists a constant $C>0$ such that the relation
$\mathbf{d}(\mathbf{x},\mathbf{y})\leq t+\varepsilon $ implies the
relation
$\mathbf{d}(\mathbf{y},\overline{\mathbf{B}}(\mathbf{x},t))\leq
C\varepsilon $ for all $t,\varepsilon >0$. It is obvious that the
strong corona property implies the corona property and that the
corona property implies the weak corona property. On the other
hand, all these properties are not equivalent to each other.

Now we can prove a theorem that relates the strong corona property
with the continuity of the map $\mathcal{B}$.

\begin{theorem}
Let $\mathbf{(X,d)}$ be a metric space and let us consider on
$\mathbf{X}\times [0,\infty)$ the metric
$d_{\infty}((\mathbf{x},t),(\mathbf{y},s))=\max\{\mathbf{d}(\mathbf{x},\mathbf{y}),|t-s|\}$.
Then, $\mathbf{(X,d)}$ has the strong corona property if and only
if $\mathcal{B}$ is a Lipschitz map.
\end{theorem}

%TCIMACRO{
%\TeXButton{margen}{\noindent%
%}}%
%BeginExpansion
\noindent%
%
%EndExpansion
\textbf{Proof. }Let $(%
\mathbf{x},t),(\mathbf{y},s)\in \mathbf{X\times [}0,\infty )$ and
let us set $\varepsilon=d_\infty((\mathbf{x},t),(\mathbf{y},s))$. Let $\mathbf{%
z\in \overline{\mathbf{B}}(x,}t\mathbf{).}$ Then,
\[
\mathbf{d(z,y})\leq t+\varepsilon \leq s+|t-s|+\varepsilon \leq
s+2\varepsilon
\]
and it follows from the strong corona property that $\mathbf{d(z,\overline{%
\mathbf{B}}(y,}s\mathbf{))}\leq 2C\varepsilon $ for certain
constant $C>0.$ Hence
\[
\sup_{\mathbf{z\in \overline{\mathbf{B}}(x,}t\mathbf{)}}\mathbf{d(z,%
\overline{\mathbf{B}}(y,}s\mathbf{))}\leq 2C\varepsilon .
\]
The same argument proves that
\[
    \sup_{\mathbf{z\in \overline{\mathbf{B}}(y,}s\mathbf{)}}
    \mathbf{d(z, \mathbf{\overline{\mathbf{B}}(x,}t\mathbf{)})}\leq
    2C\varepsilon,
\]
so that $\mathcal{B}$ is Lipschitz.

Let us now assume that $\mathcal{B}$ is Lipschitz with constant
$C$. This means that
\[
\mathbf{H}\left(\mathbf{B}(\mathbf{x},t),\mathbf{B}(\mathbf{y},s)\right)\leq
Cd_{\infty}((\mathbf{x},t),(\mathbf{y},s))
\]
for all
$(\mathbf{x},t),(\mathbf{y},s)\in\mathbf{X}\times[0,\infty)$.
Clearly,
\begin{eqnarray*}
\sup_{\mathbf{z}\in\mathbf{B}(\mathbf{x},t+\varepsilon)}\mathbf{d}(\mathbf{z},\mathbf{B_d}(\mathbf{x},t))&=&
\max\{\sup_{\mathbf{z}\in\mathbf{B}(x,t)}\mathbf{d}(\mathbf{z},\mathbf{B}(\mathbf{x},t+\varepsilon)),
\sup_{\mathbf{z}\in\mathbf{B}(\mathbf{x},t+\varepsilon)}\mathbf{d}(\mathbf{z},\mathbf{B}(\mathbf{x},t))\}\\
&=& \mathbf{H}\left(\mathbf{B}(\mathbf{x},t),\mathbf{B}(\mathbf{x},t+\varepsilon)\right)\\
&\leq & C d_{\infty}((\mathbf{x},t),(\mathbf{x},t+\varepsilon))\\
&=& C\varepsilon.
\end{eqnarray*}
Therefore,  $y\in\mathbf{B}(\mathbf{x},t+\varepsilon)$ implies
that
\[
\mathbf{d}(\mathbf{y},\mathbf{B}(\mathbf{x},t))\leq
\sup_{\mathbf{z}\in\mathbf{B}(\mathbf{x},t+\varepsilon)}\mathbf{d}(\mathbf{z},\mathbf{B}(\mathbf{x},t))\leq
C\varepsilon.
\]
This ends the proof.  $\Box $

The following examples prove that: \begin{itemize}
\item There are metric spaces which have the corona property (but not the strong
corona property) and such that $\mathcal{B}$ is not continuous.
\item There are metric spaces which do not
have the corona property and such that $\mathcal{B}$ is a
continuous map.
\end{itemize}

\begin{example}
Let us consider on $\mathbf{X}$ the discrete metric and assume that $\mathbf{%
X}$ has at least two elements. Then  $(\mathbf{X},\mathbf{d})$ has
the corona property, since for all  $\mathbf{x},\mathbf{y}\in
\mathbf{X}$ with $\mathbf{x}\neq \mathbf{y}$
and $\mathbf{d}(\mathbf{x},\mathbf{y})\leq t+\varepsilon $, we have that: i) If $t\geq 1$, then $%
\overline{\mathbf{B}}(\mathbf{x},t)=\mathbf{X}$ and $0=\mathbf{d}(\mathbf{y},\overline{\mathbf{B}}%
(\mathbf{x},t))<\varepsilon $, and ii) If $t<1$ then $t+\varepsilon <1$ for all $%
\varepsilon <\varepsilon _{0}(t)=1-t$. Hence $\mathbf{y}=\mathbf{x}$ and $0=\mathbf{d}(\mathbf{y},\overline{%
\mathbf{B}}(\mathbf{x},t))<\varepsilon $.

Now, the Hausdorff  metric $\mathbf{H}$ is the discrete metric on
$\mathbf{CLB(X)}$, so that all subsets of $\mathbf{CLB(X)}$ are
open and $\mathcal{B}$ is not continuous since $\mathcal{B}^{-1}(\{\mathbf{X}%
\})=\mathbf{X\times }[1,\infty )$ is not an open subset of $\mathbf{X\times }%
[0,\infty ).$
\end{example}

\begin{example}
Let us set $\mathbf{X}=\mathbb{R}$ with the metric given by
\[
\mathbf{d}(x,y)=\frac{|x-y|}{1+|x-y|}.
\]
Then
\[
\mathcal{B}(x,t)=[x-\frac{t}{1-t},x+\frac{t}{1-t}],
\]
if $t<1$ and $\mathcal{B}(x,t)=\mathbb{R}$ for all $t\geq 1$.
Hence
\begin{eqnarray*}
\mathbf{H}(\mathcal{B}(x,t),\mathcal{B}(y,s))&=&
\max\{\mathbf{d}(x-\frac{t}{1-t},y-\frac{s}{1-s}),\mathbf{d}(x-\frac{t}{1-t},y+\frac{s}{1-s}),\\
&\ &
\mathbf{d}(x+\frac{t}{1-t},y-\frac{s}{1-s}),\mathbf{d}(x+\frac{t}{1-t},y+\frac{s}{1-s})\},
\end{eqnarray*}
if $t,s<1$ and $\mathbf{H}(\mathcal{B}(x,t),\mathcal{B}(y,s))=1$
otherwise. Now it is easy to check that $\mathcal{B}$ is
continuous (but not uniformly continuous), since, if we fix
$(x,t)$, with $t<1$, then all the functions $\mathbf{d}(x\pm
\frac{t}{1-t},y\pm\frac{s}{1-s})$ are continuous. On the other
hand, $(\mathbb{R},\mathbf{d})$ does not satisfy the corona
property since ???

\end{example}

\section{Sufficient conditions for the weak corona property}
In this section we study several results which give us sufficient
conditions for a metric space $\mathbf{(X,d)}$ to have the weak
corona property.
\begin{theorem}
Let us assume that $(\mathbf{X},\mathbf{d})$ is boundedly compact.
Then $\mathbf{d}$ has the weak corona property.
\end{theorem}

%TCIMACRO{
%\TeXButton{margen}{\noindent%
%}}%
%BeginExpansion
\noindent%
%
%EndExpansion
\textbf{Proof }Let us assume that
$\mathbf{d}(\mathbf{x},\mathbf{y}_{n})\leq t+\varepsilon _{n}$,
for
all $n\in \mathbb{N}$, with $\varepsilon _{n}\searrow 0$, and the sequence $\{\mathbf{d}(%
\overline{\mathbf{B}}(\mathbf{x},t),\mathbf{y}_{n})\}_{n=1}^{\infty
}$ does not converge
to zero. Then there exists some $\delta >0$ and a subsequence $%
\{\mathbf{y}_{n_{k}}\}_{k=1}^{\infty }$ of
$\{\mathbf{y}_{n}\}_{n=1}^{\infty }$ such that
\[
\mathbf{d}(\overline{\mathbf{B}}(\mathbf{x},t),\mathbf{y}_{n_{k}})\geq
\delta \text{ for all }k\in \mathbb{N}\text{.}
\]
Now, $\{\mathbf{y}_{n_{k}}\}_{k=1}^{\infty }\subset \overline{\mathbf{B}}%
(\mathbf{x},t+\varepsilon _{1})$ so that there exists a subsequence $%
\{\mathbf{y}_{n_{k_{s}}}\}_{t=1}^{\infty }$ which converges to some $\mathbf{y}\in \mathbf{X}$%
, since $(\mathbf{X},\mathbf{d})$ is boundedly compact. Hence
\[
\mathbf{d}(\mathbf{x},\mathbf{y})\leq
\mathbf{d}(\mathbf{x},\mathbf{y}_{n_{k_{s}}})+\mathbf{d}(\mathbf{y}_{n_{k_{s}}},\mathbf{y})
\]
for all $s$ and, if we take the limit as $s\rightarrow \infty $, we have that $%
\mathbf{y}\in \overline{\mathbf{B}}(\mathbf{x},t)$, so that
\[
\lim_{s\rightarrow \infty }\mathbf{d}\left(\overline{\mathbf{B}}(\mathbf{x},t),\mathbf{y}_{n_{k_{s}}}\right)=0%
\text{,}
\]
a contradiction. $\Box $

\begin{corollary}
Let $\mathbf{X}$ be a finite dimensional topological manifold and
assume that $\mathbf{d}$ is a metric on $\mathbf{X}$ compatible
with its topology and $( \mathbf{X},\mathbf{d})$ is a complete
metric space. Then $\mathbf{d}$ has the weak corona property.
\end{corollary}

Being boundedly compact is a strong assumption. Can it be weakened
in some sense?. Another property which is connected to compactness
but is weaker is approximative compactness. A set
$\mathbf{A\subset X}$ is said to be
approximatively compact if for each $\mathbf{x\in X}$ and each sequence $%
\{\mathbf{a}_{n}\}_{n=1}^{\infty }\subset \mathbf{A}$ such that
\[
\lim_{n\rightarrow \infty
}\mathbf{d}(\mathbf{x},\mathbf{a}_{n})=\mathbf{d}\mathbf{(x,A):=}\inf_{\mathbf{a}\in
\mathbf{A}}\mathbf{d}(\mathbf{x},\mathbf{a})
\]
(such sequences are called minimizing sequences with respect to
$\mathbf{x}$
and $\mathbf{A}$) there exists a convergent subsequence $\{\mathbf{a}_{n_{k}}%
\}_{k=1}^{\infty }\rightarrow \mathbf{a}\in \mathbf{A}$. Now, the
following result can easily be proved:

\begin{theorem}
Let $\mathbf{(X},\mathbf{d})$ be a metric space and let us assume
that the sets
\[
\mathbf{B}(\mathbf{x},t\mathbf{)}^{c}:=\{\mathbf{y}\in
\mathbf{X}:\mathbf{d}(\mathbf{x},\mathbf{y})\geq t\}
\]
are approximatively compact for all $\mathbf{x}\in \mathbf{X}$ and
$t>0$. Then, $\mathbf{d}$ has the weak corona property.
\end{theorem}

%TCIMACRO{
%\TeXButton{margen}{\noindent%
%}}%
%BeginExpansion
\noindent%
%
%EndExpansion
\textbf{Proof }It follows from the proof of Theorem 1 that if
$\mathbf{d}$ has not the weak corona property then there exists a
sequence $\{\mathbf{y}_{n}\}_{n=1}^\infty\subset
\mathbf{B}(\mathbf{x},t)^c$ for a certain choice of $\mathbf{x}\in \mathbf{X}$ and $%
t>0$, such that $\lim_{n\rightarrow \infty }\mathbf{d}(\mathbf{x},\mathbf{y}_{n})=t$ and $\mathbf{d}(\overline{%
\mathbf{B}}(\mathbf{x},t),\mathbf{y}_{n})\geq \delta $ for all $n\in \mathbb{N}$. Now, $%
\{\mathbf{y}_{n}\}_{n=1}^{\infty }$ is a minimizing sequence with
respect to $\mathbf{x}\in
\mathbf{X}$ and $\mathbf{B}(\mathbf{x},t\mathbf{)}^{c}$, since $\mathbf{d}(\mathbf{x},\mathbf{B(}\mathbf{x},t%
\mathbf{)}^{c})=t$. Hence there exists a subsequence $\{\mathbf{y}_{n_{k}}\}_{k=1}^{%
\infty }$ which converges to some point $\mathbf{y}\in
\mathbf{B(}\mathbf{x},t\mathbf{)}^{c}$.
Clearly, $\mathbf{y}\in \overline{\mathbf{B}}(\mathbf{x},t)$ and $\mathbf{d}(\mathbf{y}_{n_{k}},\mathbf{y})%
\rightarrow 0$ as $k\rightarrow \infty $, a contradiction. $\Box $

\section{Topological vector spaces without the weak corona
property}

We have already proved that some topological spaces have not the
weak corona property, but these spaces had some holes.  Thus, it
is not clear, by the moment, if there exists some example of a
complete metric topological vector space (which is a nice space,
with no holes) without the weak corona property. Surprisingly, the
answer is affirmative, as the following result proves.

\begin{theorem} Let $\mathbf{X}=\{\{a_{n}\}_{n=1}^{\infty }\subseteq\mathbb{R}:\#\{k:a_{k}\neq
0\}<\infty \}$ be the vector space of all sequences of real
numbers with finite support and let us define
\[
\mathbf{d}(\{a_{n}\}_{n=1}^\infty,\{b_{n}\}_{n=1}^\infty\mathbf{)=}\sum_{n=1}^{\infty
}\Phi _{n}(a_{n}-b_{n})\text{,}
\]
where
\[
\Phi _{k}(t)=\left\{
\begin{array}{lll}
k|t| &  & \text{if }|t|\leq \frac{1}{k} \\
1+|t|-1/k &  & \text{otherwise}
\end{array}
\right.
\]
Then $\mathbf{(X,d)}$ is a metric topological vector space without
the corona property. Consequently, its metric completion $(%
\widetilde{\mathbf{X}},\mathbf{d)}$ is a complete metric
topological vector space without the corona property.
\end{theorem}
\noindent \textbf{Proof. } It is an easy exercise to prove that
$\mathbf{d}$ is a metric on $\mathbf{X}$ which is invariant by
translations and that $(\mathbf{X,d)}$ is a topological vector
space (i.e., sum and scalar product are continuous). On the other
hand,
\[
\mathbf{d}\left(\{0\}_{n=1}^\infty,\{{\textstyle\frac{2}{k}}\delta
_{k,n}\}_{n=1}^{\infty }\right)=1 +\frac{1}{k}
\]
and
\[
\mathbf{d}\left(\{{\textstyle\frac{2}{k}}\delta
_{k,n}\}_{n=1}^{\infty
},\mathbf{\overline{B}}(\{0\}_{n=1}^\infty,1)\right)=1
\]
for all $k$.

To prove this claim, let $\{a_{n}\}_{n=1}^\infty\in
\overline{\mathbf{B}}(\{0\}_{n=1}^\infty,1 )$ be arbitrarily
chosen. Then $\sum_{n=1}^{\infty }\Phi _{n}(a_{n})\leq 1$. This
implies that $\Phi _{n}(a_{n})\leq 1$ for all $n$. Hence
$|a_{n}|\leq \frac{1}{n}$ for all $n$ and $|\frac{2}{k}-a_{k}|\geq
\frac{1}{k}$, so that
\[
\mathbf{d}\left(\{0\}_{n=1}^\infty,\{{\textstyle\frac{2}{k}}\delta
_{k,n}\}_{n=1}^\infty \right)\geq \Phi _{k}(\frac{2}{k}\delta
_{k,k}-a_{k})\geq 1.
\]
On the other hand, $\{\frac{1}{k}\delta _{k,n}\}_{n=1}^{\infty }\in \overline{\mathbf{B}}%
(\{0\}_{n=1}^\infty,1)$ and $\mathbf{d}\left(\{\frac{2}{k}\delta
_{k,n}\}_{n=1}^{\infty },\{\frac{1}{k}\delta
_{k,n}\}_{n=1}^{\infty }\right)= 1$. Hence $(\mathbf{X,d)}$ does
not satisfy the weak corona property. This ends the proof. $\Box$

\section{Function spaces, the weak corona property and an open question}

 The space of continuous functions $\mathbf{X=}C(\mathbb{%
R})$ is a metrizable t.v.s. In fact, if we denote by
$\{p_{n}\}_{n=0}^{\infty }$ the family of seminorms
$p_{n}(f)=\max_{x\in \lbrack -n,n]}|f(x)|$, $n\in \mathbb{N}$, and we
take $\{c_{n}\}_{n=0}^\infty\searrow 0^{+}$ an arbitrary
nonincreasing sequence of real numbers which converges to zero,
then  $\mathbf{d}(f,g)=\sup_{n\in \mathbb{N}}\{c_{n}%
\frac{p_{n}(f-g)}{1+p_{n}(f-g)}\}$ is a metric on $C(\mathbb{R})$
which is invariant by translations. Now, we are able to prove the
following:

\begin{theorem}
$(C(\mathbb{R}),\mathbf{d})$ has the weak corona property.
\end{theorem}

%TCIMACRO{
%\TeXButton{margen}{\noindent%
%}}%
%BeginExpansion
\noindent%
%
%EndExpansion
\textbf{Proof }It follows from the invariance by translations of the metric,
that we can assume (without loss of generality) that our balls are centered
in the null function. Now, $\mathbf{d}(f,\mathbf{0})\leq t$ if and only if
\[
\frac{p_{n}(f)}{1+p_{n}(f)}\leq \frac{t}{c_{n}}\text{ for all }n\in
\{k:c_{k}>t\}
\]
which is equivalent to say that
\[
p_{n}(f)\leq \frac{t}{c_{n}-t}\text{ for all }n\leq n_{0}(t);
\]
where $\{k:c_{k}>t\}=\{1,2,...,n_{0}(t)\}$.

Let $t>0$ be arbitrarily fixed and let $\varepsilon
_{0}=\varepsilon_{0}(t)>0$ be such that
\[
c_{k}>t\Leftrightarrow c_{k}>t+\varepsilon
\]
for all $\varepsilon \leq \varepsilon _{0}$. With this choice of
$\varepsilon_{0}$, we have that
$f_{\varepsilon}\in\mathbf{\overline{B}(0},t+\varepsilon)$
$(\varepsilon\leq \varepsilon_{0})$ means that the graph of
$f_{\varepsilon}$ is contained in certain boxes
$\mathbf{B}_{k,\varepsilon}=[-k,k]\times J_{k,\varepsilon}$
$(k\leq n_0(t))$ which are a little bit taller than the boxes
$\mathbf{B}_{k}=[-k,k]\times J_{k}$ $(k\leq n_0(t))$
 which contain the graph of an arbitrary function
 $g\in \mathbf{\overline{B}}(0,t)$ (see the picture).

%\begin{picture}(200,270)(-15,65)
\begin{center}
    \includegraphics[scale=0.4]{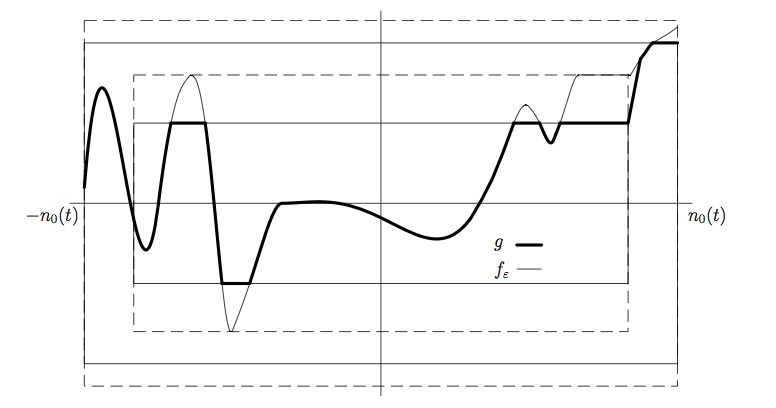}
    \end{center}
    %\put(-128,173){$g$}
    %\put(-128,155){$f_\varepsilon$}
    %\put(-430,190){$-n_0(t)$}
    %\put(-3,190){$n_0(t)$}
%    \put(-75,115){$t$}
%    \put(-85,60){$d(y_n,\mathbf{B}_d(x,t))$}
%    \put(-135,96){$\mathbf{B}_d(x,t)$}
%    \put(-105,85){$x$}
%    \put(-25,85){$y_n$}
%\end{picture}

%\begin{center}
%\includegraphics[scale=0.5]{dibujo1.eps}
%\end{center}
Let us assume that $\mathbf{d}%
(f_{\varepsilon },\mathbf{0})\leq t+\varepsilon $, and $\varepsilon \leq
\varepsilon _{0}/2$. Then it is easy to find a continuous function $g\in
\mathbf{C[}-n_{0}(t),n_{0}(t)\mathbf{]}$ such that
\begin{eqnarray*}
p_{n}(f_{\varepsilon }-g) &\leq &\max_{k\leq n}\{\frac{t+\varepsilon }{%
c_{k}-t-\varepsilon }-\frac{t}{c_{k}-t}\} \\
&=&\max_{k\leq n}\frac{c_{k}\varepsilon }{(c_{k}-t-\varepsilon )(c_{k}-t)}%
\text{ for all }n\leq n_{0}(t)=n_{0}(t+\varepsilon )
\end{eqnarray*}
This function can be extended with continuity to all the real line
satisfying that $h(x)=|(f_{\varepsilon }-g)(x)|$ is a nonincreasing function on the intervals $[n_0(t),n_0(t)+1/2]$ and $[-n_0(t)-1/2,-n_0(t)]$, and vanishes for $|x|>n_{0}(t)+1/2$. Then we obtain that
\[
p_{n_{0}(\varepsilon )+j}(f_{\varepsilon }-g)=p_{n_{0}(\varepsilon
)}(f_{\varepsilon }-g)\leq M(t,\varepsilon ):=\max_{k\leq n_{0}(t)}\frac{%
c_{k}\varepsilon }{(c_{k}-t-\varepsilon )(c_{k}-t)}\text{ for all }j\in \mathbb{%
N}
\]
and
\[
c_{n}\frac{p_{n}(f_{\varepsilon }-g)}{1+p_{n}(f_{\varepsilon }-g)}\leq c_{n}%
\frac{M(t,\varepsilon )}{1+M(t,\varepsilon )}\leq c_{0}\frac{M(t,\varepsilon
)}{1+M(t,\varepsilon )}\text{ for all }n\in \mathbb{N}
\]
This obviously implies that
\[
\mathbf{d}(f_{\varepsilon },\mathbf{\mathbf{\mathbf{\overline{\mathbf{B}}%
}}(0,}t\mathbf{)})\leq \mathbf{d}(f_{\varepsilon },g)\leq c_{0}\frac{%
M(t,\varepsilon )}{1+M(t,\varepsilon )}\rightarrow 0\text{ (for }\varepsilon
\rightarrow 0\text{)}
\]
and the space $C(\mathbb{R})$ satisfies the weak corona property. $\Box $

 It is not difficult to check that all
classical topological function vector spaces satisfy the weak corona
property, because of analogous arguments to those given in the above proof.
On the other hand, the situation in more general contexts is not clear by the moment. This lead us to the following

%TCIMACRO{
%\TeXButton{margen}{\noindent%
%}}%
%BeginExpansion
\noindent%
%
%EndExpansion
\textbf{Open question: }Let $\mathbf{X}$ be a locally convex
topological vector space and assume that its topology is the
associated to the separating family of seminorms
$\{p_{n}\}_{n=0}^{\infty }$. Let
$\mathcal{C}=\{c_{n}\}_{n=0}^\infty$ be a sequence of positive
real numbers such that $\mathcal{C}$ is non-increasing,
$\lim_{n\rightarrow \infty }c_{n}=0$. We define the metric
\[
\mathbf{d}(x,y)=\sup_{n\in \mathbb{N}}\{c_{n}\frac{p_{n}(x-y)}{%
1+p_{n}(x-y)}\}\text{ }
\]
Is it true that $(\mathbf{X,d})$ satisfies the weak corona property?

\end{document}